\documentclass[12pt]{amsart}

\usepackage{mathrsfs}
\usepackage{amsmath,amssymb,amsfonts,latexsym,txfonts}
\usepackage{eucal}
\usepackage{fancyhdr}

\newlength{\originalbase}
\newtheorem{theorem}{Theorem}[section]
\newtheorem{proposition}[theorem]{Proposition}

\newtheorem{definition}[theorem]{Definition}
\newtheorem{corollary}[theorem]{Corollary}

\input amssym.def
 
\input amssym.tex

\renewcommand{\SS}{\mathbb{S}}
\newcommand{\RR}{\mathbb{R}}

\newcommand{\LL}{\mathcal{L}}
\newcommand{\KK}{\mathscr{K}}

\begin{document}

\begin{center}
{\large\bf If you can hide behind it, can you hide inside it?}\\[12mm]
\end{center}

\begin{center}
{\bf Daniel A. Klain} \\
Department of Mathematical Sciences\\
University of Massachusetts Lowell\\
Lowell, MA 01854 USA  \\
Daniel\_{}Klain@uml.edu\\
\end{center}

\begin{quotation}
{\small
\noindent
{\bf Abstract \/}  
Let $K$ and $L$ be compact convex sets in $\RR^n$.
Suppose that, for a given dimension $1 \leq d \leq n-1$, 
every $d$-dimensional orthogonal projection of $L$
contains a translate of the corresponding projection of $K$.
Does it follow that the original set $L$ contains a translate of $K$?  
In other words, if $K$ can be translated to ``hide behind" $L$ from any perspective,
does it follow that $K$ can ``hide inside" $L$?

A compact convex set $L$ is defined to be $d$-{\em decomposable} if $L$ is a direct Minkowski sum (affine Cartesian
product) of two or more convex bodies each of dimension at most $d$.  
A compact convex set $L$ is called $d$-{\em reliable} if, whenever each $d$-dimensional orthogonal projection of $L$ 
contains a translate of the corresponding $d$-dimensional projection of $K$, 
it must follow that $L$ contains a translate of $K$.

It is shown that, for $1 \leq d \leq n-1$:
\begin{enumerate}
\item[(1)]  $d$-decomposability implies $d$-reliability.  

\item[(2)] A compact convex set $L$ in $\RR^n$ is $d$-reliable if and only if,
for all $m \geq d+2$,
no $m$ unit normals to regular boundary points of $L$ form the
outer unit normals of a $(m-1)$-dimensional simplex. 

\item[(3)] Smooth convex bodies are not $d$-reliable.

\item[(4)] A compact convex set $L$ in $\RR^n$ is $1$-reliable if and only if $L$ is 
$1$-decomposable (i.e.~a parallelotope).  \sloppy

\item[(5)] A centrally symmetric compact convex set $L$ in $\RR^n$ is $2$-reliable if and only if $L$ is 
$2$-decomposable. \\ 
However, there are non-centered $2$-reliable convex bodies that
are not $2$-decomposable. 
\end{enumerate}
As a result of (5) above, 
the only reliable centrally symmetric covers in $\RR^3$ from the perspective of 2-dimensional
shadows are the affine convex cylinders (prisms).  
However, in dimensions greater than 3, it is shown
that 3-decomposability is only sufficient, and not necessary, for $L$ to cover reliably 
with respect to $3$-shadows,
even when $L$ is assumed to be centrally symmetric.
}
\end{quotation}
\vspace{6mm}

Consider two compact convex subsets $K$ and $L$ of 
$n$-dimensional Euclidean space.  
Suppose that, for a given dimension $1 \leq d < n$,
every $d$-dimensional orthogonal projection (shadow) of $L$
contains a translate of the corresponding projection of $K$.
Does it follow that the original set $L$ contains a translate of $K$?  
In other words, if $K$ can be translated to ``hide behind" $L$ from any perspective,
does it follow that $K$ can ``hide inside" $L$?

In dimension 2 it is easy to see that the answer is {\em No}. 
For example, if an equilateral triangle $\Delta$ is inscribed in a disc $D$ of unit diameter, the 
slightly larger triangle $(1+\epsilon) \Delta$
still has less than unit width in every direction (provided $\epsilon > 0$ is sufficiently small), but
no longer fits inside $D$.  The same construction works for any set $K$ inscribed in $D$ and having
strictly less than unit diameter.  Another counterexample arises from comparing $\Delta$
and the dilated and reflected triangle $-(1+\epsilon)\Delta$ for small $\epsilon > 0$.

Although the details are less obvious, counterexamples also exist in
higher dimensions.  
Let $B$ denote the unit Euclidean 3-ball, 
and let $T$ denote the regular tetrahedron having edge length $\sqrt{3}$.  
Jung's Theorem \cite[p. 84]{Bonn2}\cite[p. 320]{Webster}
implies that every 2-projection of $T$ 
is covered by a translate of the unit disk.  But a simple computation shows that 
$B$ cannot cover the tetrahedron $T$.  
An analogous construction yields a similar result for higher dimensional simplices and Euclidean balls. 
One might say that, although $T$ can be translated within a fixed distance from $B$ 
(i.e. without moving far away) to hide behind $B$ from any observer's perspective, 
this does not imply that $T$ can hide {\em inside} $B$.  

Indeed, for $1 \leq d \leq n-1$,
it is shown in \cite{Klain-Inscr} that if $K$ is a compact convex set in $\RR^n$
having at least $d+2$ exposed points, then there exists another compact convex set $L$
such every $d$-dimensional orthogonal projection (shadow) of $L$
contains a translate of the corresponding projection of $K$,
while $L$ does not contain a translate of $K$.
In certain cases one can even find examples
where $K$ also has {\em larger volume} than $L$ (and so certainly could not fit inside $L$).
For a detailed example of this volume phenomenon, see \cite{Klain-Shadow}.

This leads to the question: under what additional conditions on 
the covered set $K$ or the covering set $L$ does covering of shadows of $K$
by shadows of $L$
imply covering of the original set $K$ by the set $L$? 

This question is easily answered when a sufficient degree of symmetry is imposed.
For example, a support function argument implies that the answer is {\em Yes}
if {\em both} of the bodies $K$ and $L$ are centrally symmetric.
It is also not difficult to show that if every $d$-projection of $K$ (for some $1 \leq d < n$)
can be translated into the corresponding shadow of an orthogonal $n$-dimensional box $C$,
then $K$ fits inside $C$ by some translation, since one needs only to check that the widths
are compatible in the $n$ edge directions of $C$.  A similar observation applies if $C$
is a parallelotope (an affine image of a box), or even a cylinder (the product of
an $(n-1)$-dimensional compact convex set with a line segment).

In \cite{Lut-contain} Lutwak uses Helly's theorem to
prove that, if every $n$-simplex containing $L$ also contains a translate of $K$, then
$L$ contains a translate of $K$.  In the present article we
generalize Lutwak's theorem
in order to reduce questions about shadow covering to questions about 
circumscribing simplices and simplicial cylinders.  
A compact convex set $L$ will be called $d$-{\em decomposable} if $L$ is a {\em direct} Minkowski sum (affine Cartesian
product) of two or more convex bodies each of dimension at most $d$ (see Section~\ref{csas}).  
A compact convex set $L$ will be called $d$-{\em reliable} if, whenever each $d$-shadow of $L$ 
contains a translate of the corresponding $d$-shadow of $K$, 
it follows that $L$ contains a translate of $K$
(see Section~\ref{when}).  It will be shown that, for $1 \leq d \leq n-1$:
\begin{enumerate}
\item[(1)]  $d$-decomposability implies $d$-reliability.  
(Theorem~\ref{lutshadow2})

\item[(2)] A compact convex set $L$ in $\RR^n$ is $d$-reliable if and only if,
for all $m \geq d+2$,
no $m$ unit normals to regular boundary points of $L$ form the
outer unit normals of a $(m-1)$-dimensional simplex. (Theorem~\ref{rel})

\item[(3)] Smooth convex bodies are {\em not} $d$-reliable.
(Corollary~\ref{smooth})

\item[(4)] A compact convex set $L$ is $1$-reliable if and only if $L$ is 
$1$-decomposable (i.e.~a parallelotope).  (Corollary~\ref{1iff1})

\item[(5)] A {\em centrally symmetric} compact convex set $L$ is $2$-reliable if and only if $L$ is 
$2$-decomposable. (Theorem~\ref{2iff2})\\ 
However, there are non-centered $2$-reliable convex bodies that
are not $2$-decomposable.  (Corollary~\ref{polyrely})
\end{enumerate}
As a result of (5) above, 
the only reliable centrally symmetric covers in $\RR^3$ from the perspective of 2-dimensional
shadows are the affine convex cylinders (prisms).  
However, in dimensions greater than 3, it will be seen
(at the end of Section~\ref{cscs})
that 3-decomposability is only sufficient, and not necessary, for $L$ to cover reliably 
with respect to $3$-shadows,
even when $L$ is assumed to be centrally symmetric.

The containment and covering problems addressed in this article are special cases of
the following general question: under what conditions 
will a compact convex set necessarily contain a translate or otherwise congruent copy of another?  
Progress on different aspects of this general question also appears in the work of
Gardner and Vol\v{c}i\v{c} \cite{Gard-Vol},
Groemer \cite{Gro-proj},
Hadwiger \cite{Had3,Had4,Had-proj,Lincee,Santa}, 
Jung \cite{Bonn2,Webster}, 
Lutwak \cite{Lut-contain}, 
Rogers \cite{Rogers}, 
Soltan \cite{Soltan},
Steinhagen \cite[p. 86]{Bonn2}, 
Zhou \cite{Zhou1,Zhou2}, and many others (see also \cite{Gard2006}).

\section{Background} 
Denote $n$-dimensional Euclidean space by $\RR^n$, and let $\SS^{n-1}$ denote the set of
unit vectors in $\RR^n$; that is, the unit $(n-1)$-sphere centered at the origin.

Denote by $\KK_n$ the set of compact convex subsets of $\RR^n$.  The $n$-dimensional
(Euclidean) volume of a convex set $K$ will be denoted $V_n(K)$.  If $u$ is a unit vector in
$\RR^n$, denote by $K_u$ the orthogonal projection of a set $K$ onto the subspace $u^\perp$.
More generally, if $\xi$ is a $d$-dimensional subspace of $\RR^n$, denote by $K_\xi$
the orthogonal projection of a set $K$ onto the subspace $\xi$.  The boundary of a compact
convex set $K$ will be denoted by $\partial K$.

Let $h_K: \RR^n \rightarrow \RR$ denote the support function of a compact convex set $K$;
that is,
$$h_K(v) = \max_{x \in K} x \cdot v.$$
If $\xi$ is a subspace of $\RR^n$ then the support function $h_{K_\xi}$ is given by the restriction
of $h_K$ to $\xi$.
If $u$ is a unit vector in
$\RR^n$, denote by $K^u$ the support set of $K$ in the direction of $u$; that is,
$$K^u = \{x \in K \; | \; x \cdot u = h_K(u) \}.$$
If $P$ is a convex polytope, then $P^u$ is the face of $P$ having $u$ in its outer normal cone.

Given two $K, L \in \KK_n$ and $a,b \geq 0$ denote
$$aK + bL = \{ax + by \; | \; x \in K \hbox{ and } y \in L\}.$$
An expression of this form is called a {\em Minkowski combination} or 
{\em Minkowski sum}.  Because $K$ and $L$ are convex, the set $aK + bL$ is also convex.  
Convexity also implies that $aK + bK = (a+b)K$ for all $a,b \geq 0$.
Support functions satisfy the identity
$h_{aK+bL} = ah_K + bh_L$.  (See, for example, any of 
\cite{Bonn2,red,Webster}).

If $K \in \KK_n$ has non-empty interior, define the {\em surface area measure} $S_K$ on
the $(n-1)$-dimensional unit sphere $\SS^{n-1}$ as follows:  For $A \subseteq \SS^{n-1}$
denote by $K^A = \bigcup_{u \in A} K^u$, and define
$S_K(A) = \mathcal{H}_{n-1}(K^A)$, the $(n-1)$-dimensional Hausdorff measure of the
subset $K^A$ of the boundary of $K$. (See \cite[p. 203]{red}.)
If $P$ is a polytope, then $S_P$ is a pointed measure concentrated 
at precisely those directions $u$ that are outer normals to the facets of $P$.

The measure $S_K$ is easily shown to satisfy the property 
\begin{equation}
\int_{\SS^{n-1}} u \; dS_K = \vec{o},
\label{surf}
\end{equation}
that is, the mass distribution on the sphere described by $S_{K}$ has center of mass
at the origin.  For a convex polytope $P$ having outward facet unit normals $u_1, \ldots, u_m$
and corresponding facet areas $\alpha_1, \ldots, \alpha_m > 0$, the identity~(\ref{surf}) takes
the simple and intuitive form:
\begin{equation}
\alpha_1 u_1 + \cdots + \alpha_m u_m = \vec{o}.
\label{surfp}
\end{equation}

Minkowski's Existence Theorem \cite[p. 125]{Bonn2}\cite[p. 390]{red} gives a 
useful converse to the identity~(\ref{surf}):  If $\mu$ is a non-negative measure on the unit sphere $\SS^{n-1}$ such that $\mu$ has center of mass at the origin, and if $\mu$ is not concentrated on any great (equatorial) $(n-1)$-subsphere, then $\mu = S_K$ for some $K \in \KK_n$.  Moreover, this convex body $K$ is {\em unique up to translation}.


Suppose that $\mathscr{F}$ is a family of compact convex sets in $\RR^n$.
Helly's Theorem \cite{Bonn2,red,Webster} asserts that if every $n+1$ sets
in $\mathscr{F}$ share a common point, then the entire family shares a common point.
In \cite{Lut-contain} Lutwak used Helly's theorem to
prove the following fundamental criterion for whether a set $L \in \KK_n$ contains 
a translate of another set $K \in \KK_n$.   
\begin{theorem}[Lutwak's containment theorem] 
Let $K,L \in \KK^n$.  
The following are equivalent:
\begin{enumerate}
\item[\bf (i)] For every simplex $\Delta$ such that $L \subseteq \Delta,$ 
there exists $v \in \RR^n$ such that 
$K + v \subseteq \Delta$.
\item[\bf (ii)] There exists $v_0 \in \RR^n$ such that
$K + v_0 \subseteq L$.
\end{enumerate}
\label{lut}
\end{theorem}
In other words, if every $n$-simplex containing $L$ also contains a translate of $K$, then
$L$ contains a translate of $K$.

\section{Circumscribing sets and shadows}
\label{csas}

A convex set $K \in \KK_n$ will be called 
{\em $d$-decomposable} if there exists a subspace decomposition
$$\RR^n = \xi_1 \oplus \cdots \oplus \xi_m,$$
where $\dim \xi_i \leq d$ for each $i$, 
and compact convex sets $K_i \subseteq \xi_i$ for each $i$,
such that  $K = K_1 + \cdots + K_m$.  Decompositions of this kind will be denoted
$$K = K_1 \oplus \cdots \oplus K_m.$$
If $K = \Delta_1 \oplus \cdots \oplus \Delta_m$, where the component sets $\Delta_i$ are {\em simplices},
each of dimension at most $d$,
then we will say that $K$ is a {\em $d$-decomposable simplex product}.
The product will be called {\em orthogonal} if the subspaces $\xi_i$ are mutually
orthogonal.

The $2$-decomposable sets in $\RR^3$ (as well as products of $(n-1)$-dimensional sets with line
segments in $\RR^n$) are often called {\em cylinders} or {\em prisms}.  

If the circumscribing simplices are replaced 
by circumscribing simplex products for $L$, then the following generalization of
Lutwak's Theorem~\ref{lut} is obtained.
\begin{theorem}[Prismatic containment theorem] 
Let $K,L \in \KK^n$.  
The following are equivalent:
\begin{enumerate}
\item[\bf (i)] For every $d$-decomposable set $C \in \KK_n$ such that $L \subseteq C,$ 
there exists $v \in \RR^n$ such that 
$K + v \subseteq C$.
\item[\bf (ii)] For every $d$-decomposable orthogonal simplex product $C$ such that $L \subseteq C,$ 
there exists $v \in \RR^n$ such that 
$K + v \subseteq C$.
\item[\bf (iii)] For every $d$-dimensional subspace $\xi \subseteq \RR^n$, there exists $w \in \xi$
such that $K_\xi + w \subseteq L_\xi$.
\end{enumerate}
\label{lutshadow}
\end{theorem}
In other words, if every $d$-decomposable (simplex) product $C$ 
containing $L$ also contains a translate of $K$, then
every $d$-shadow $L_\xi$ contains a translate of the corresponding shadow $K_\xi$, and vice versa.  

The following proposition will simplify the proof of Theorem~\ref{lutshadow}.
\begin{proposition} Let $K, L \in \KK_n$.  Let $\psi: \RR^n \rightarrow \RR^n$ be a
nonsingular linear transformation.  Then $L_u$ contains a translate of $K_u$ for all unit directions $u$
if and only if $(\psi L)_u$ contains a translate of $(\psi K)_u$ for all $u$.
\label{obl}
\end{proposition}
This proposition implies that nothing is gained (or lost) by allowing
more general (possibly non-orthogonal) linear projections.

\begin{proof} 
For $S \subseteq \RR^n$ and a nonzero vector $u$, let $\LL_S(u)$ denote the set of straight lines in $\RR^n$
parallel to $u$ and meeting the set $S$.
The projection $L_u$ contains a translate $K_u$ for each unit vector $u$ if and only if,
for each $u$, there exists $v_u$ such that
\begin{equation}
\LL_{K+v_u} (u) \subseteq \LL_L(u).
\label{hum}
\end{equation}
But $\LL_{K+v_u}(u) = \LL_{K}(u) + v_u$ 
and $\psi \LL_{K}(u) = \LL_{\psi K}(\psi u)$.
It follows that~(\ref{hum}) holds if and only if
$\LL_{K}(u) + v_u \subseteq \LL_L(u)$,
if and only if
$$\LL_{\psi K}( \psi u) + \psi v_u \subseteq \LL_{\psi L}(\psi u) \;\;\; \hbox{ for all unit } u,$$
Set 
$$\tilde{u} = \frac{\psi u}{|\psi u|} \;\;\; \hbox{ and } \;\;\; \tilde{v} = \psi v_u.$$
The relation~(\ref{hum}) now holds if and only if,
for all $\tilde{u}$, there exists $\tilde{v}$ such that
$$\LL_{\psi K}(\tilde{u}) + \tilde{v} \subseteq \LL_{\psi L}(\tilde{u}),$$
which holds if and only if 
$(\psi L)_{\tilde{u}}$ contains a translate of $(\psi K)_{\tilde{u}}$ for all $\tilde{u}$.
\end{proof}

\begin{proof}[{\bf\em Proof of Theorem~\ref{lutshadow}}] 
To begin, note that {\bf (i)} implies {\bf (ii)} trivially.

Suppose that {\bf (ii)} holds.  Given a $d$-subspace $\xi \subseteq \RR^n$, let $T$ be a simplex
in $\xi$ that circumscribes $L_\xi$.  Let $u_{d+1}, \ldots, u_n$ be an orthonormal basis for $\xi^\perp$, and 
let $C'$ be a cube in $\xi^\perp$ with edges parallel to the directions $u_i$ and large enough so that
$L_{\xi^\perp} \subseteq C'$.  Now let $C = T \oplus C'$.
Since $L \subseteq C$, it follows from {\bf (ii)} that there exists $v \in \RR^n$ such that
$K+v \subseteq C$.  This implies that $K_\xi + v_\xi \subseteq C_\xi = T$.  On applying Lutwak's Theorem~\ref{lut}
in the subspace $\xi$ it follows that $K_\xi + w \subseteq L_\xi$ for some $w \in \xi$.  Therefore, {\bf (ii)}
implies {\bf (iii)}.

Next, suppose that {\bf (iii)} holds.   If $L \subseteq C = C_1 \oplus \cdots \oplus C_m$,
where each $\dim \xi_i \leq d$, then let $\psi$ be a non-singular linear operator on $\RR^n$ such 
that the subspaces $\psi(\xi_i)$ are mutually orthogonal.  By Proposition~\ref{obl}, the condition
{\bf (iii)} also holds for $\psi K$ and $\psi L$.  For each $i$ we obtain $v_i \in \psi \xi_i$ such
that 
$$(\psi K)_{\psi \xi_i} + v_i \subseteq (\psi L)_{\psi \xi_i} \subseteq 
(\psi C)_{\psi \xi_i} = \psi C_i.$$  
Let $v = v_1 + \cdots + v_m$.  
Since the subspaces $\psi\xi_i$ are mutually orthogonal, we have
$\psi K+v \subseteq \psi C$, so that $K + \psi^{-1}v \subseteq C$.  
Therefore, {\bf (iii)} implies {\bf (i)}, and the three assertions are equivalent.
\end{proof}

It is worth noting the following special case of Theorem~\ref{lutshadow}.
\begin{theorem} 
Let $K,C \in \KK^n$, where $C$ is $d$-decomposable.

Suppose that, for each $d$-dimensional subspace $\xi \subseteq \RR^n$, 
there exists $w \in \xi$
such that $K_\xi + w \subseteq C_\xi$.
Then there exists $v \in \RR^n$ such that $K+v \subseteq C$.
\label{lutshadow2}
\end{theorem}
When $d = n-1$, Theorem~\ref{lutshadow2} says that if you can hide behind a cylinder from
any perspective (and without rotating), then you can also hide inside the cylinder.

More consequences of Theorem~\ref{lutshadow} are explored in \cite{Klain-Shadow}.  

\section{Simplicial families of unit normals}

Theorem~\ref{lutshadow2} motivates a converse question: 
If $L$ is {\em not} $d$-decomposable, does there 
necessarily exist $K$ such that every $d$-shadow of $L$ contains a translate of the
corresponding $d$-shadow of $K$, while $L$ itself does not contain a translate of $K$?  
The answer is {\em not necessarily.}  
We will show in a later section (see Corollary~\ref{polyrely})
that if $L$ is a square pyramid (the convex hull of
a square in $\RR^3$ with a point above its center) then 
no $K$ can hide behind $L$ unless $K$ can also hide inside $L$.  
However, the square pyramid is not $2$-decomposable.
In other words, the condition of being $d$-decomposable is sufficient, but not necessary.

In this section we develop some tools for constructing necessary and sufficient conditions for
when shadow covering implies actual covering.  These tools are applied in later sections.

A set of unit vectors $\{ u_0, \ldots, u_d \}  \subseteq \SS^{n-1}$ 
will be called a {\em $d$-simplicial family},
or {\em $d$-simplicial}, if $u_0, \ldots, u_d$ 
span a $d$-dimensional subspace of $\RR^n$, and if there exist real numbers
$c_0, \ldots, c_d > 0$, such that
$$c_0 u_0 + c_1 u_1 + \cdots + c_d u_d = o.$$
Equivalently, $u_0, \ldots, u_d$ are the outer unit normals of some $d$-dimensional simplex.
Note that a $d$-simplicial family contains exactly $d+1$ unit vectors.

It will be seen in Sections~\ref{when} and~\ref{cscs} that
certain translative covering properties of a compact convex set $L$ hinge in the
existence of simplicial families of unit normals to regular points of $L$.  
The next three propositions will be used in that context.
(Readers in a hurry may wish to scan Sections~\ref{when} and~\ref{cscs} 
and return to these technical points later on.)

\begin{proposition} Suppose that $A = \{u_1, \ldots, u_m\} \subseteq \SS^{n-1}$ 
contains no simplicial families of size 3 or greater, and that
\begin{equation}
c_1 u_1 + \cdots + c_m u_m = o
\label{orvex}
\end{equation}
for some $c_1, \ldots, c_m > 0$.  Then $m = 2s$ for some integer $s$, and 
there exist linearly independent vectors $v_1, \ldots, v_s \in \SS^{n-1}$, where $s \leq n$, such
that
$$A = \{\pm v_1, \ldots, \pm v_s\}.$$
\label{one}
\end{proposition}

\begin{proof}  By~(\ref{orvex}) the set $A$ must have at least 2 elements, and 
if $A$ has size $2$ then the proposition is trivial.  

Suppose that the proposition fails for some set $A$ of minimal size $m$, 
where $m > 2$.  By~(\ref{orvex}) there 
exists a minimal subfamily $\{ u_{i_1}, \ldots, u_{i_k} \} \subseteq A$ such that
$$a_1 u_{i_1} + \cdots + a_k u_{i_k} = o$$
for some $a_1, \ldots, a_k > 0$.  Let $d = \dim(\mathrm{Span}\{u_{i_1}, \ldots, u_{i_k}\})$.  Since the
$u_{i_j}$ are dependent, we have $k \geq d+1$.   
If $k > d+1$ then Caratheodory's theorem \cite[p. 3]{red}
(applied in the span of the $\{u_{i_j}\}$) implies that the origin $o$ lies in the convex hull of a sub-subfamily of size at most $d+1 < k$ of the
$u_{i_j}$, violating the minimality of $k$.  Therefore $k = d+1$, and $\{ u_{i_1}, \ldots, u_{i_k} \}$ is a simplicial set.  By the original assumption on simplicial families in $A$, it follows that $k=2$,
so that $u_i = -u_j$ for some $i \neq j$.  

Without loss of generality, suppose that $u_1 = - u_2$ and that $c_1 \geq c_2$. 
It now follows from ~(\ref{orvex}) that
$$o = (c_1 - c_2) u_1 + c_3 u_3 + \cdots + c_m u_m$$
Suppose $c_1 - c_2 > 0$.  The minimality of $m$ implies that the proposition holds for the set
$\{ u_1, u_3, \ldots, u_m \}$, so that $m-1$ is even and these remaining vectors $u_1, u_3, \ldots, u_m$ 
can be partitioned into
distinct antipodal pairs.  Since $u_1 = -u_2$, this would violate the original assumption that the $u_i$
are distinct.  Therefore $c_1 = c_2$, and
$$o = c_3 u_3 + \cdots + c_m u_m.$$
Once again the minimality of $m$ implies that
that the proposition holds for the set
$u_3, \ldots, u_m$, so that $m-2$ is even (and therefore $m$ is even), and the remaining $u_i$ 
can be separated into distinct antipodal pairs $\pm v_2, \ldots, \pm v_s$, 
where the $v_i$ are linearly independent.

It remains to show that $u_1$ (and similarly $u_2$) is linearly independent from the vectors $v_i$.
If $u_1$ lies in the span of $v_2, \ldots, v_s$, then $u_1$ lies in the span of a minimal subset
of $v_2, \ldots, v_k$ of size at least 2, since $u_1$ is distinct from each $\pm v_i$.  The resulting
linear dependence relation violates the nonexistence of simplicial subsets of size 3 or greater inside $A$.  

Setting $v_1 = u_1$ now completes the proof of the proposition.
\end{proof}

\begin{proposition} Suppose that $A \subseteq \SS^{n-1}$ 
contains no simplicial sets of size 3 or greater, and that the origin lies in
the interior of the convex hull of $A$.
Then there exist linearly independent vectors $v_1, \ldots, v_n \in \SS^{n-1}$, such that
$$A = \{\pm v_1, \ldots, \pm v_n\}$$
\label{oneone}
\end{proposition}

\begin{proof}  By Caratheodory's theorem there exists a
finite subfamily $u_1, \ldots, u_m$ of $A$, such that~(\ref{orvex}) holds.  
Since $o$ lies in the {\em interior} of the convex hull of $A$, we can take $m$ large enough
so that $u_1, \ldots u_m$ spans $\RR^n$.
By Proposition~\ref{one}, this subfamily
has the form $\{\pm v_1, \ldots, \pm v_s\} \subseteq A$, where $v_1, \ldots, v_s$ are linearly independent
and span $\RR^n$.  It follows that $s=n$.

If $w \in A$ and $w \neq \pm v_i$, then
$w$ lies in the span of some $v_{i_1}, \ldots, v_{i_k}$, where $k \geq 2$ is minimal.  This 
linear dependence relation violates the nonexistence of simplicial sets of size 3 or greater inside the set $A$.  

It follows that $\{\pm v_1, \ldots, \pm v_n\} = A$.
\end{proof}

\begin{proposition}  Suppose that $A \subseteq \SS^{n-1}$ is symmetric under reflection through the
origin; that is $A = -A$.  Suppose also that $A$
contains no simplicial sets of size 4 or greater, and that the origin lies in
the interior of the convex hull of $A$.

Then there exists a subspace direct sum decomposition
$$\RR^n = W_1 \oplus \cdots \oplus W_k,$$ 
where each $\dim W_i \leq 2$, and such that
$A \subseteq W_1 \cup \cdots \cup W_k$. 
\label{two}
\end{proposition}

\begin{proof}  Since $A = -A$, the set $A$ is composed of antipodal pairs $\pm v$ of unit vectors.
Moreover, since the convex hull of $A$ has interior and is centrally symmetric, 
there exist at least $n$ pairs  
$\pm u_1, \ldots, \pm u_n$ in $A$ whose $n$ directions are linearly independent.
If $A = \{ \pm u_1, \ldots, \pm u_n \}$ then $\RR^n$ is a direct sum of the lines
spanned by each $\pm u_i$, and the proposition follows.

If, instead, $\pm v$ is another antipodal pair in $A$, where $v \neq \pm u_i$ for all $i$, then
without loss of generality (relabeling the signs on $\pm u_i$ as needed), we have
$$-v = c_1 u_1 + \cdots + c_d u_d$$
for some $c_1, \ldots, c_d > 0$, where $d$ is minimal.  If $d \geq 3$ then the relation
$$v + c_1 u_1 + \cdots + c_d u_d = 0$$
implies that $\{v, u_1, \ldots, u_d\}$ form a simplicial family in $A$
of size at least 4, contradicting hypothesis.  Meanwhile, since $v \neq \pm u_i$, we must
have $d > 1$.   The remaining possibility is $d=2$, so that $v$ lies in the span of $\{u_1, u_2\}$.

If $w \in A$ and $w \neq \pm v, \pm u_1, \ldots, \pm u_n$, 
then $w$ lies in the span of 2 of the $u_i$ by a similar argument.  
But if $w = a_1 u_1 + a_3 u_3,$ say, where $a_1, a_3 > 0$, then
$$w = \frac{a_1}{c_1} (-c_2 u_2 - v) + a_3 u_3 = 
\frac{a_1 c_2}{c_1}(-u_2) + \frac{a_1}{c_1}(-v) + a_3 u_3,$$
Since every 3 of the 4 vectors $v, w, u_2, u_3$ are linearly independent, we obtain
a simplicial set of size 4, another contradiction.  Therefore, either $w$ also lies 
in the span of $\{u_1, u_2\}$ or in the span of $\{u_i, u_j\}$ for $j>i>2$.  An iteration of this argument 
implies that $\RR^n$ is decomposed into a direct sum 
$\RR^n = W_1 \oplus \cdots \oplus W_{\lfloor \frac{n+1}{2} \rfloor}$ of subspaces $W_i$ 
each having dimension at most 2, and where every $v \in A$ also lies in some $W_i$.
\end{proof}

We will also need the following proposition, which clears up ambiguities regarding when 
shadows cover inside a larger ambient space.
\begin{proposition} Suppose that $\xi$ is a linear flat in $\RR^n$.
Let $K$ and $L$ be compact convex sets in $\xi$.  
Suppose that, for each $d$-subspace $\eta \subseteq \xi$, the projection
$L_\eta$ contains a translate of $K_\eta$.  Then
$L_\eta$ contains a translate of $K_\eta$ for every $d$-subspace $\eta \subseteq \RR^n$.
\label{include}
\end{proposition}
In other words, embedding $K$ and $L$ in a higher-dimensional space does not change whether or not every $d$-shadow of $L$ contains a translate of the corresponding $d$-shadow of $K$ (even though there are now more shadow directions to verify).

\begin{proof}  
Suppose that $\eta$ is a $d$-subspace of $\RR^n$.  Let $\hat{\eta}$ denote the orthogonal projection
of $\eta$ into $\xi$.  Since $\dim(\hat{\eta}) \leq \dim(\eta) = d$, we can
translate $K$ and $L$ inside $\xi$ so that $K_{\hat{\eta}} \subseteq L_{\hat{\eta}}$.  Let us assume this
translation has taken place.  Note that, for $v \in \hat{\eta}$, we now have
$h_K(v) \leq h_L(v)$.

If $u \in \eta$, then express $u = u_{\xi} + u_{\xi^\perp}$.  Since $K \subseteq \xi$,
$$h_K(u) = \max_{x \in K} x \cdot u 
=\max_{x \in K} x \cdot u_\xi = h_K(u_\xi),
$$
and similarly for $L$.  But since $u \in \eta$, we have $u_\xi \in \hat{\eta}$, so that
$$h_K(u) = h_K(u_\xi)
\leq h_L(u_\xi) = h_L(u).
$$
In other words, $K_\eta \subseteq L_\eta$.
\end{proof}

\section{When can a convex set conceal without covering?}
\label{when}

We now address the possibility of a converse to Theorem~\ref{lutshadow2}.
\begin{definition}  Suppose that $1 \leq d \leq n-1$.
A compact convex set $L$ in $\RR^n$ is said to be a {\em $d$-reliable cover}, 
or {\em $d$-reliable}, if whenever $K \in \KK_n$ and 
every $d$-shadow $L_\xi$ contains a translate of the corresponding shadow $K_\xi$, 
it follows that $L$ contains a translate of $K$.
\end{definition}
Evidently, if $L$ is $d$-reliable, then $L$ is also $m$-reliable for all $m > d$.

Theorem~\ref{lutshadow2} asserts that 
if $L$ is $d$-decomposable then $L$ is also $d$-reliable.
However, we will see that a square pyramid gives a counterexample to the converse assertion.  
It is $2$-reliable, but not $2$-decomposable (Corollary~\ref{polyrely}).

The next two theorems describe a
{\em necessary and sufficient} condition for $L$ to be a $d$-reliable cover.
Recall that a point $x$ on the boundary of a compact convex set $L$ is said 
to be {\em regular} if the outward normal cone
to $L$ at $x$ contains exactly one unit vector. 

\begin{theorem}  Suppose that $L$ has regular boundary points $x_0, \ldots, x_{d+1}$,
whose corresponding unit normals $u_0, \ldots, u_{d+1}$ are a simplicial family.
Then there exists a polytope $S$ such that $L_\xi$ contains a translate of $S_\xi$ for each
$d$-subspace $\xi$ of $\RR^n$, while $L$ {\em does not} contain a translate of $S$.  In particular,
$L$ is not $d$-reliable.
\label{rel1}
\end{theorem}

\begin{proof} First, note that, by Proposition~\ref{include}, it sufficient to prove this theorem for the case in which $L$ has interior.  For if $L$ lacks interior, we simply restrict our attention to the affine hull of $L$.   Once the theorem is verified in this case, one can
apply Proposition~\ref{include} to verify the theorem when $L$ is re-embedded in a higher-dimensional space.  So let us now assume that $L$ has interior.

Suppose that $L$ has regular boundary points $x_0, \ldots, x_{d+1},$ as in the hypothesis of the theorem. 
Let $S$ be the convex hull of $\{x_0, \ldots, x_{d+1}\}$.  Evidently $S \subseteq L$.  Since 
$\{u_0, \ldots, u_{d+1} \}$ is a simplicial family, there exist $c_i > 0$ such that
\begin{equation}
c_0 u_0 + \cdots + c_{d+1} u_{d+1} = o.
\label{z}
\end{equation}
Moreover, every $d+1$ of the $u_i$ are linearly independent, so that no subfamily of
the $u_i$ contains the origin in its convex hull, and the property~(\ref{z}) does not hold
for any subfamily.

Since $S \subseteq L$, we have $h_S \leq h_L$.  Moreover, by our choices of $x_i$ and $u_i$, 
$h_L(u_i) = u_i \cdot x_i \leq h_S(u_i)$.  Therefore, $h_L(u_i) = h_S(u_i)$ for each $i$.

Let $\xi$ be a $d$-dimensional
subspace of $\RR^n$.  If $\pi_\xi(x_i)$ lies on the boundary of $L_\xi$, then 
$h_{L_\xi}(w) = x_i \cdot w$ for some unit $w \in \xi$.  
Since $h_{L_\xi}$ is given by the restriction of $h_L$
to the subspace $\xi$, it follows that $h_L(w) = x_i \cdot w$.  
By the regularity of the boundary point
$x_i$, we have $w = u_i$, so that $u_i \in \xi$.

Similarly, if $u_i \in \xi$ then
$$h_{L_\xi}(u_i) =  h_{L}(u_i) = x_i \cdot u_i = \pi_\xi(x_i) \cdot u_i,$$
so that $\pi_\xi(x_i)$ lies on the boundary of $L_\xi$, with outward unit normal $u_i$.

Since the $u_i$ form a simplicial family of size $d+2$, at most $d$ of these vectors
$u_i$ can lie in $\xi$.  
It follows that $S_\xi$ meets the boundary of $L_\xi$ at $j+1$ points, for some $j \leq d-1$,
having outward normals $u_0, \ldots, u_{j}$ (without loss of generality).
Since $o$ does not lie in the convex hull of $u_0, \ldots, u_j$, 
there exists a unit vector $v \in \xi$ such
that $v \cdot u_i < 0$ for $i \leq j$.  

Let $y_i = \pi_\xi(x_i)$.
If $i \leq j$ then $y_i$ is a regular point of the boundary $\partial L_\xi$
with outward unit normal $u_i$.  
Since each $v \cdot u_i < 0$ in this case, there exists $\epsilon > 0$ such that
each $y_i + \epsilon v$ lies in the relative interior of $L_\xi$.

Meanwhile, if $i > j$ then $y_i$ lies in the relative interior of $L_\xi$ already, so that
$y_i + \epsilon v$ lies in the relative interior of $L_\xi$ as well,
provided we have chosen $\epsilon > 0$ small enough.
In other words, there exists $\epsilon > 0$ so that
$y_i + \epsilon v$ lies in the relative interior of $L_\xi$ for all $i$.

Let $T = S_\xi + \epsilon v$. 
Since the polytope $T$ is the convex hull of the
points $y_i + \epsilon v$,
it follows that $T$ lies in the relative interior of $L_\xi$.  Therefore,
there exists $a_\xi > 1$ such that $L_\xi$ contains a translate of $a_\xi T$,
whence $L_\xi$ contains a translate of $a_\xi S_\xi$.

Since the set of all $d$-subspaces of $\RR^n$ is compact, there exists $\alpha > 1$, 
independent of $\xi$, such
that some translate of $\alpha S_\xi$ lies inside $L_\xi$ for each $\xi$.  

On the other hand, if $\alpha S+w \subseteq L$ for some $w$, then
$$h_L(u_i) \geq h_{\alpha S+w}(u_i) = \alpha h_S(u_i) + u_i \cdot w = 
\alpha h_L(u_i) + u_i \cdot w > h_L(u_i) + u_i \cdot w,$$
so that $u_i \cdot w < 0$ for all $i$.  This strict inequality contradicts~(\ref{z}).
\end{proof}

To prove the converse to Theorem~\ref{rel1}, we first consider the polytope case.
Recall that a {\em facet} of a polytope $Q$ is a face of co-dimension 1 in the affine
hull of $Q$.
\begin{theorem}
Let $K \in \KK_n$, and let $Q$ be a convex polytope in $\RR^n$.  Suppose that
$Q_\xi$ contains a translate of $K_\xi$ for every $d$-subspace $\xi$, and that
$Q$ does not contain a translate of $K$.
Then there exists a simplicial family of facet unit normals $\{ u_0, \ldots, u_m \}$ 
to $Q$, for some $m \geq d+1$.
\label{rel2}
\end{theorem}
In other words, if a polytope $Q$ is {\em not} $d$-reliable, then $Q$ has a simplicial
family of facet unit normals of size at least $d+2$.

\begin{proof} As in the previous proof, 
Proposition~\ref{include} makes it sufficient to verify the case 
in which $Q$ has interior.  

Suppose that $K_\xi$ can be translated inside $Q_\xi$ for each
$d$-subspace $\xi$, while $K$ cannot be translated inside $Q$.
Without loss of generality, translate $K$ so that the origin $o$ lies inside
the relative interior of $K$.  This implies that $h_K \geq 0$.
  
Since $Q$ has interior,
there exists $\epsilon > 0$ such that $\epsilon K$ can be translated inside $Q$.
Since $Q$ is compact we may assume $\epsilon$ to be maximal.  Evidently $\epsilon < 1$,
since no translate of $K$ fits inside $Q$.  
Without loss of generality, translate $Q$ so that $\epsilon K \subseteq Q$.

Denote the facets of $Q$ by $F_0, \ldots, F_q$, having outward unit normals
$u_0, \ldots, u_q$.  Suppose that $\epsilon K$ meets facets $F_0, \ldots, F_m$, and
misses the others.

If the convex hull of $\{u_0, \ldots u_m\}$ does not contain the origin $o$, then 
there exists a vector $v$ such that $v \cdot u_i < 0$ for $i = 0, \ldots, m$.  This implies
that, for suffiiciently small $\delta$, the translate $\epsilon K + \delta v$ lies in the
{\em interior} of $Q$.  This violates the maximality of $\epsilon$.   Therefore, there exist
$a_0, \ldots, a_m \geq 0$ such that
$$ a_0 u_0 + \cdots + a_m u_m = o. $$
Renumbering the facets as necessary, we have
\begin{equation}
c_0 u_0 + \cdots + c_s u_s = o
\label{zz}
\end{equation}
where each $c_i > 0$ and $s$ is minimal, so that $\{u_0, \ldots, u_s\}$ is a simplicial family.

If $s \leq d$, then the $s+1$ unit vectors $u_i$ lie inside a $d$-subspace $\xi$.
Since $\epsilon K$ meets each of the facets $F_0, \ldots, F_s$, we have
\begin{equation}
\epsilon h_K(u_i) = h_{\epsilon K}(u_i) = h_Q(u_i) = h_{Q_{\xi}}(u_i)
\label{allsame}
\end{equation}
for each $i = 0, \ldots, s$.  Since $Q_\xi$ contains a translate of $K_\xi$, there exists $w \in \xi$
so that $K_\xi + w \subseteq Q_\xi$, and
$$\epsilon h_K(u_i) = h_{Q_\xi}(u_i) \geq h_{K_\xi}(u_i) + w \cdot u_i 
= h_K(u_i) + w \cdot u_i$$
for each $i = 0, \ldots, s$.
After summing over $i$, it follows from~(\ref{zz}) that
$$\epsilon \sum_{i=0}^s c_i h_K(u_i) \geq 
\sum_{i=0}^s c_i h_K(u_i) + w \cdot \sum_{i=0}^s c_i u_i
=  \sum_{i=0}^s c_i h_K(u_i).$$
Recall that $h_K \geq 0$ and each $c_i > 0$.
Since $\epsilon < 1$, it follows that
$$\sum_{i=0}^s c_i h_K(u_i) = 0,$$
so that each $h_K(u_i) = 0$.  Therefore, each
$h_Q(u_i) = 0$, by~(\ref{allsame}).  It now follows from~(\ref{zz}) and the 
sublinearity of the support function $h_Q$ that
the projection of $Q$ onto the span of $\{ u_0, \ldots, u_s \}$ is a single point.
This is 
a contradiction, since $Q$ has interior.
It follows that $s \geq d+1$.

Therefore, there exists a simplicial family of facet unit normals $u_0, \ldots, u_s$ to $Q$, 
where $s \geq d+1$.
\end{proof}

Putting Theorems~\ref{rel1} and~\ref{rel2} together, we obtain the following.
\begin{theorem}[Reliability Theorem] Let $L \in \KK_n$.  Then $L$ is a $d$-reliable cover 
if and only if every simplicial family of normals 
to regular boundary points of $L$ has size at most $d+1$.
\label{rel}
\end{theorem}

\begin{proof} 
Suppose a simplicial family of unit normals 
to regular boundary points of $L$ has size $d+2$ or greater.   By Theorem~\ref{rel1}, $L$ is not $d$-reliable.

To prove the converse, suppose that $L$ is not $d$-reliable.  
Then there exists $K \in \KK_n$ such that $L_\xi$ contains a translate of $K_\xi$
for every $d$-subspace $\xi$, while $L$ does not contain a translate of $K$.

Since regular points are dense on the boundary of $L$ (see \cite[p. 73]{red}), 
there exists a countable dense set of regular
points on the boundary of $L$.  By intersecting half-spaces that support $L$ at these points, construct 
a sequence of polytopes $P_i$, decreasing with respect to set inclusion,
such that $P_i \rightarrow L$ and each $P_i$ has facet normals that are 
unit normals at regular points of $L$.

If $P_i$ contains a translate of $K$ for all $i$, then so does $L$, a contradiction.  Therefore,
there exists $j$ such that $P_j$ does not contain a translate of $K$.  But each projection
$L_\xi \subseteq (P_j)_\xi$, so that each projection $(P_j)_\xi$ contains a translate of $K_\xi$.
In other words, the polytope $P_j$ is not $d$-reliable.  By Theorem~\ref{rel2}, there are facet
unit normals $u_0, \ldots, u_m$ for the polytope $P_j$ that form a simplicial family, for some
$m \geq d+1$.  Since the facet normals of $P_j$ were taken from unit normals to regular
points of $L$, this completes the proof.
\end{proof}

Recall that a simplex $T$ {\em circumscribes} $L$ if $L \subseteq T$ 
and if $aT$ contains no translate of $L$ when $a < 1$.  An $n$-simplex $T \supseteq L$ circumscribes
$L$ if and only if $L$ meets every facet of $T$.  Theorem~\ref{rel} therefore implies the following. 
\begin{corollary} Let $L \in \KK_n$.  Then $L$ is $(n-1)$-reliable if and only if there is no circumscribing
$n$-simplex $T$ of $L$ such that $\partial T \cap L$ consists of regular points.
\end{corollary}

Since every boundary point of a {\em smooth} convex body is a regular boundary point,
the following corollary is now immediate.
\begin{corollary}  If $L$ is a smooth convex body in $\RR^n$,
there exists an $n$-simplex $S$ such that $L_u$ contains a translate of $S_u$ for every
unit direction $u$, while $L$ {\em does not} contain a translate of $S$.
\label{smooth}
\end{corollary}

We can now characterize $1$-reliability.
\begin{corollary} A convex set $L \in \KK_n$ is a $1$-reliable cover if and only if $L$ is a parallelotope.
\label{1iff1}
\end{corollary}

\begin{proof}  If $L$ is a parallelotope then $L$ is $1$-reliable, by Theorem~\ref{lutshadow2}.

Conversely, if $L$ is $1$-reliable, then Theorem~\ref{rel} asserts that 
there are no simplicial sets of size $3$ or more among the unit normals at regular points of $L$.  
By Proposition~\ref{include} we may assume, without loss of generality, that $L$ has interior.
In this case there exist affinely independent unit normals 
$u_1, \ldots, u_{m}$ at regular points of $L$, 
where $m \geq n+1$, and
where the $u_i$ 
do not all lie in the same hemisphere.  It follows that
$$o = c_1 u_1 + \cdots + c_{m} u_{m}$$
for some $c_1, \ldots, c_{m} > 0$.  By Proposition~\ref{oneone}, the set of regular normals
of $L$ has the form $\{\pm v_1, \ldots, \pm v_n\}$, for some
linearly independent set $v_1, \ldots, v_n \in \SS^{n-1}$.  
Let $P$ be the unique (up to translation)
parallelotope having facet unit normals $\pm v_i$ and corresponding facet areas $c_i$.  
Since the regular points of $L$ are dense in the boundary of $L$, it follows that
$L$ and $P$ must be translates.
\end{proof}

\begin{corollary} Suppose that $P$ is a polytope in $\RR^n$.  Then $P$ is a $d$-reliable cover
if and only if, for all $m \geq d+2$,
no $m+2$ facets of $P$ share normal directions with an $(m-1)$-simplex.
\label{polyrely}
\end{corollary}

Since no four facet normals of the square pyramid $P$ in $\RR^3$ contain the origin in the interior of
their convex hull, any $K \in \KK_3$ that can ``hide behind" $P$ can also ``hide inside" $P$.
In other worlds, $P$ is 2-reliable, in spite of being indecomposable.

\section{Centrally symmetric covering sets}
\label{cscs}

We saw in the previous section that $L$ is a reliable 1-cover if and only if $L$ is 
1-decomposable (i.e.~a  parallelotope).
However, the square pyramid is 2-reliable in spite of being indecomposable.  

A compact convex set $L$ is said to be {\em centrally symmetric} if $L$ and $-L$ are translates.
For equivalence 
of 2-reliability and 2-decomposability to hold, we must restrict our attention to centrally symmetric bodies.

\begin{theorem}  A centrally symmetric set $L \in \KK_n$ is $2$-reliable if and only if $L$ is  2-decomposable.
\label{2iff2}
\end{theorem}

The 3-dimensional case of Theorem~\ref{2iff2} has the following especially simple form.
\begin{corollary}  A centrally symmetric set $L \in \KK_3$ is $2$-reliable if and only if $L$ is a cylinder.
\end{corollary}

The proof of Theorem~\ref{2iff2} will use the following auxiliary results.

\begin{proposition} Let $P$ be a convex polytope in $\RR^n$ with non-empty interior.
Suppose that $\xi$ is a proper subspace of $\RR^n$, 
and suppose that each facet unit normal of $P$ lies either in $\xi$ or in $\xi^\perp$.

Then there exist 
polytopes $P_1 \subseteq \xi$ and $P_2 \subseteq \xi^\perp$ such that
$P = P_1 \oplus P_2$.
\label{o-exist}
\end{proposition}

\begin{proof}  
Suppose that the facet unit normals of $P$ are given by 
$$\{u_1, \ldots, u_p, v_1, \ldots, v_q\}$$
where $u_1, \ldots, u_p \in \xi$ and 
$v_1, \ldots v_q \in \xi^\perp$.  Suppose that each facet of $P$ with normal $u_i$ has area $a_i$ and
each facet with normal $v_j$ has area $b_j$.  
By the Minkowski condition,
$$a_1 u_1 + \cdots + a_p u_p + b_1 v_1 + \cdots + b_q v_q = o$$
It follows from the independence of $\xi$ and $\xi^\perp$ that
$$a_1 u_1 + \cdots + a_p u_p = o \quad \hbox{ and } \quad b_1 v_1 + \cdots + b_q v_q = o$$

By the Minkowski Existence Theorem \cite{Bonn,red} there exists a polytope $Q_1 \subseteq \xi$ having
facet normals $u_i$ and corresponding facet areas $a_i$.
Similarly, there exists a polytope $Q_2 \subseteq \xi^\perp$ having
facet normals $v_j$ and corresponding facet areas $b_j$.  

Let $d = \dim \xi$, so that $\dim \xi^\perp = n-d$.
For $x,y > 0$, the Minkowski sum $xQ_1 + yQ_2$ has the same unit normals as $P$ and has corresponding
facets $xQ_1^{u_i} + yQ_2$
and $xQ_1 + yQ_2^{v_j}$, having the respective facet areas $x^{d-1} y^{n-d}  V_{n-d}(Q_2) a_i$
and $x^d y^{n-d-1}  V_d(Q_1) b_j$.
Set
$$x = \left( \frac{V_{n-d}(Q_2)^{n-d-1}}{V_{d}(Q_1)^{n-d}} \right)^{\frac{1}{n-1}}
\quad \hbox{ and } \quad
y = \left( \frac{V_{d}(Q_1)^{d-1}}{V_{n-d}(Q_2)^{d}} \right)^{\frac{1}{n-1}},$$
and let $P_1 = xQ_1$ and $P_2 = y Q_2$.  The polytope $P_1 \oplus P_2$ now has the same facet normals
and the same corresponding facet areas as $P$.  It follows from the uniqueness assertion of
the Minkowski Existence Theorem that $P$ and $P_1 \oplus P_2$ must be translates.
\end{proof}

\begin{proposition} Let $K \in \KK_n$ have non-empty interior.  Suppose that
there is a subspace decomposition $\RR^n = \xi \oplus \xi'$ 
such that each unit normal at a regular point of $K$ 
lies either in $\xi$ or in $\xi'$.

Then there is a subspace decomposition $\RR^n = \eta \oplus \eta'$, where $\dim \eta = \dim \xi$
and $\dim \eta' = \dim \xi'$, and compact 
convex sets $K_1 \subseteq \eta$ and $K_2 \subseteq \eta'$ such that
$K = K_1 \oplus K_2$.
\label{exist2}
\end{proposition}

\begin{proof}  To begin, suppose that $\xi' = \xi^\perp$, so that $\RR^n = \xi \oplus \xi'$ is
an orthogonal decomposition.
Since regular points are dense on the boundary of $K$ (see \cite[p. 73]{red}), 
there exists a countable dense set of regular
points on the boundary of $K$.  By intersecting half-spaces that support $K$ at these points, construct 
a sequence of polytopes $P_i$, decreasing with respect to set inclusion,
such that $P_i \rightarrow K$ and each $P_i$ has facet normals that are 
unit normals at regular points of $K$.  

By Proposition~\ref{o-exist}, each $P_i = Q_i \oplus Q'_i$, where $Q_i \subseteq \xi$ and $Q'_i \subseteq \xi'$.
Since projections are continuous, the $Q_i = (P_i)_\xi$ converge to $K_\xi$, and 
similarly $Q'_i \rightarrow K_{\xi'}$.
Therefore $K = \lim_{i} P_i = K_\xi \oplus K_{\xi'}$.

More generally, if $\xi$ and $\xi'$ are not orthogonal complements, then
let $\psi: \RR^n \rightarrow \RR^n$ 
be a nonsingular linear transformation such that $\psi^{-T} \xi \perp \psi^{-T} \xi'$,
where  $\psi^{-T}$ denotes 
the inverse transpose of $\psi$.  
Let $\eta = \psi^{-1} \psi^{-T} \xi$ and $\eta' = \psi^{-1} \psi^{-T} \xi'$.

Recall that $(\psi \xi)^\perp = \psi^{-T} (\xi^\perp)$.
Therefore, if each
unit normal at a regular point of $K$ 
lies either in $\xi$ or in $\xi'$, then each
unit normal at a regular point of $\psi K$ 
lies either in $\psi^{-T} \xi$ or in $\psi^{-T} \xi'$.  Since these subspaces form an orthogonal
decomposition,
the previous argument implies that $\psi K = L_1 \oplus L_2$,
where $L_1 \subseteq \psi^{-T} \xi$ and $L_2 \subseteq \psi^{-T} \xi'$. 
It follows that $K = K_1 \oplus K_2$, where $K_1 = \psi^{-1} L_1 \subseteq \eta$ and 
$K_2 = \psi^{-1} L_2 \subseteq \eta'$.  
\end{proof}

\begin{proof}[{\bf\em Proof of Theorem~\ref{2iff2}}]
If $L$ is $2$-decomposable then $L$ is $2$-reliable by Theorem~\ref{lutshadow2}.

For the converse, suppose that $L$ is $2$-reliable.
Let $A$ denote the set of unit normals at regular points of $L$.
Since $L$ is 2-reliable, $A$ contains no simplicial subsets sets of size $4$, 
by Theorem~\ref{rel}.  

Since $L$ is centrally symmetric, we have $A = -A$.  By Proposition~\ref{two}, 
there exists a subspace direct sum decomposition
$$\RR^n = W_1 \oplus \cdots \oplus W_k,$$ 
where each $\dim W_i \leq 2$, and such that
$A \subseteq W_1 \cup \cdots \cup W_k$. 
It follows from Proposition~\ref{exist2} that $L$ is $2$-decomposable.
\end{proof}

In view of Theorem~\ref{2iff2}, one may be tempted to conjecture that $d$-reliability 
is equivalent to $d$-decomposability for centrally symmetric bodies, but this turns out to be false for $d=3$.  
Consider the following 12 vectors in $\RR^4$:
$$\pm(1, 1, 0, 0)
\pm(1, 0, 1, 0), \;
\pm(1, 0, 0, 1), \;
\pm(0, 1, 1, 0), \;
\pm(0, 1, 0, 1), \;
\pm(0, 0, 1, 1)$$
By Minkowski's existence theorem, there exists a unique 12-faceted polytope $Q$ in $\RR^4$, centrally symmetric about the origin (i.e. $Q = -Q$), having facet normals parallel to the directions above, with each facet
having unit 3-volume.   One can verify that the set of vectors above contains no simplicial 5-family,
so that $Q$ is 3-reliable by Corollary~\ref{polyrely}.  
A routine linear algebra computation (using Proposition~\ref{exist2}) 
also verifies that $Q$ is not $3$-decomposable.

\section{Some open questions}

There remain
several fundamental open questions about convex bodies and projections, among them the following:
\begin{enumerate}
\item[I.]  Under what symmetry (or other) 
conditions on $L \in \KK_n$ is $d$-reliability
equivalent to $d$-decomposability, for $d > 2$?
\end{enumerate}
A solution to Problem I would generalize Corollary~\ref{1iff1} and Theorem~\ref{2iff2}.  For example, what happens
if we assume that $L$ is a zonoid?

Denote the $n$-dimensional (Euclidean) volume of $L \in \KK_n$ by $V_n(L)$.  
\begin{enumerate}
\item[II.]  
Let $K, L \in \KK_n$ such that $V_n(L) > 0$, and let $1 \leq d \leq n-1$.
Suppose that $L_\xi$ contains a translate of $K_\xi$ 
for every $d$-subspace $\xi$ of $\RR^n$.  \\

\noindent
What is the best upper bound for the ratio $\frac{V_n(K)}{V_n(L)}$?
\end{enumerate}
Some partial answers to Problem II are offered in \cite{Klain-Shadow}.  There it is shown that if $K_\xi$ can be translated inside $L_\xi$ for all $d$-dimensional subspaces $\xi$, then 
then $K$ has smaller volume than $L$ whenever $L$ can be approximated by Blaschke combinations of
$d$-decomposable sets.  However, there are cases in which $V_n(K) > V_n(L)$, in spite of the covering
condition on shadows.
For more $L \in \KK_n$ it is also shown that, if $K_u$ can be translated inside $L_u$ for all unit directions $u$, then $V_n(K) \leq nV_n(L)$, where $n$ is the dimension of the ambient space for $K$ and $L$.  However, I doubt this is the best possible bound.  

\begin{enumerate}
\item[III.] Let $K, L \in \KK_n$, and let $1 \leq d \leq n-1$.
Suppose that, for each $d$-subspace $\xi$ of $\RR^n$, the orthogonal projection 
$K_\xi$ of $K$ can be moved inside $L_\xi$ by some {\em rigid motion}
(i.e.~a combination of translations, rotations, and reflections).  \\

\noindent
Under what simple (easy to state, easy to verify) additional conditions
does it follow that $K$ can be moved inside $L$ by a rigid motion?
\end{enumerate}
Problem III 
is an intuitive generalization of the questions addressed in this
article.  Indeed, each question can be re-phrased allowing for rotations (and reflections) as well as
translations.  
However, the arguments 
presented so far rely on the observation that the set of
translates of $K$ that fit inside $L$, that is, the set
$$\{v \in \RR^n \; | \; K+v \subseteq L\},$$
is a compact convex set in $\RR^n$.  
By contrast, the set of rigid motions
of $K$ that fit inside $L$ will lie in a more complicated Lie group.  For this reason (at least)
the questions of covering via rigid motions may be more difficult to address than
the case in which only translation is allowed.

\bibliographystyle{amsplain}
\providecommand{\bysame}{\leavevmode\hbox to3em{\hrulefill}\thinspace}
\providecommand{\MR}{\relax\ifhmode\unskip\space\fi MR }
\providecommand{\MRhref}[2]{%
  \href{http://www.ams.org/mathscinet-getitem?mr=#1}{#2}
}
\providecommand{\href}[2]{#2}

\end{document}